\documentclass[12pt]{article}
\usepackage{amssymb}

\def\appendix{\par}  
\def\sn{{\mathfrak s \mathfrak n}}
\def\zz{{\mathfrak z}}
\def\cc{{\mathfrak c}}
\def\sqbdd{{\sqsubseteq^{{\rm bounded}}}}

\def\Sig{{\bf \Sigma_1^1}}
\def\Sigtwo{{\bf \Sigma_2^0}}
\def\Sigthree{{\bf \Sigma_3^0}}
\def\Pig{{\bf \Pi_1^1}}
\def\Pigtwo{{\bf \Pi_2^0}}
\def\Pigthree{{\bf \Pi_3^0}}
\def\trees{{\mathcal T}}
\def\poset{{\mathbb P}}
\def\pos{{\mathbb P}(B)}
\def\lq{\mbox{``}}
\def\rq{\mbox{''}}
\def\sq{\subseteq}
\def\bb{{\mathfrak b}}
\def\zz{{\mathfrak z}}
\def\om{\omega}
\def\infsub{[\om]^\om}
\def\finsub{[\om]^{<\om}}
\def\proof{\par\noindent Proof\par\noindent}
\def\name#1{\stackrel{\circ}{#1}}

\def\forces{{| \kern -2pt \vdash}}

\def\res{\upharpoonright}
\def\qed{\par\noindent QED\par}
\def\rmand{\mbox{ and }}
\def\rmor{\mbox{ or }}

\newtheorem{theorem}{Theorem}
\newtheorem{lemma}[theorem]{Lemma}

\newtheorem{prop}[theorem]{Proposition}

\begin{document}

\begin{center}
{\large On relatively analytic and Borel subsets}
\end{center}

\begin{flushright}
Arnold W. Miller\footnote{
Thanks to University of Florida, Gainesville and to Boise State University,
Idaho for their hospitality during the time this paper
was written and to J.Zapletal and T.Bartoszynski for some helpful discussions. 
}
\end{flushright}

\begin{center}
Abstract
\end{center}

\begin{quote}

Define $\zz$ to be the smallest cardinality of a function $f:X\to Y$ with
$X,Y\sq 2^\om$ such that there is no Borel function $g\supseteq f$.  In this
paper we prove that it is relatively consistent with ZFC to have $\bb<\zz$
where $\bb$ is, as usual, smallest cardinality of an unbounded family in
$\om^\om$.  This answers a question raised by Zapletal.

We also show that it is relatively consistent with ZFC that there exists
$X\sq 2^\om$ such that the Borel order of $X$ is bounded but there
exists a relatively analytic subset of $X$ which is not relatively coanalytic.
This answers a question of Mauldin.
\end{quote}

The following is an equivalent definition of $\zz$:
$$\zz=\min\{|X|:X\sq 2^\om, \exists Y\sq X\;\;\; Y \mbox{ is not Borel in }X\}$$
For one direction 
we can use for each $Y\subseteq X$ its characteristic function $f:X\to 2$.  
For the other direction use that a function is Borel iff the inverse
image of each basic open set is Borel.

The following answers a question raised by Zapletal \cite{zap} see
appendix A.

\begin{theorem}\label{zapletal}
It is relatively consistent with ZFC that $\bb<\zz$.
\end{theorem}

Define $p\in \poset(A)$ for $A\sq 2^\om$ iff 
$p$ is a finite set of consistent sentences of
the form: 
\begin{enumerate}
\item $\lq x\in \cap_{m<\om}U_{nm}\rq $ where $x\in A$, $n\in\om$, or
\item $\lq x\notin U_{nm}\rq $ where $x\in 2^\om$, $n,m\in\om$, or
\item $\lq [s]\sq U_{nm}\rq $ where $s\in 2^{<\om}$, $n,m\in\om$.
\end{enumerate} 
By consistent we simply mean the following: 
\begin{itemize}
\item $p$ cannot contain both $\lq x\in \cap_{m<\om}U_{nm}\rq$ and $\lq x\notin U_{nk}\rq$ 
 for some  $x,n,k$, and
\item $p$ cannot contain both $\lq x\notin U_{nm}\rq$ and 
$\lq [x\res k]\sq U_{nm}\rq$
for some $x.n,m,k$.
\end{itemize} 
The ordering on $\poset(A)$ is given by inclusion: $p\leq q$ iff $p\supseteq q$.
Note that the set $A$ enters into the picture only in sentence of 
type (1).  

This partial order is from Miller \cite{borhier} where there are 
versions for all
countable Borel orders (this is for $\Sigthree$). It can be looked
on as a generalization of almost disjoint forcing of Jensen and Solovay. I
learned about describing almost disjoint forcing as sets of sentences from 
Jack Silver.

Now suppose that $G$ is $\poset(A)$-generic over $V$.  Define
$$U^G_{nm}=\cup\{[s]: \lq [s]\sq U_{nm}\rq\in G\}
\rmand W^G_n=\cap_{m<\om}U_{nm}^G$$

\begin{lemma}\label{one}
 For any $x\in V\cap 2^\om$
\begin{enumerate}
\item $x\notin U^G_{nm}$ iff $\lq x\notin U_{nm}\rq\in G$
\item $x\in W_n^G$ iff $\lq x\in\cap_{m<\om} U_{nm}\rq\in G$
\item $x\in A$ iff $x\in \cup_{n<\om}W_n^G$
\end{enumerate}
\end{lemma}
\proof
To prove $(1)$ working in $V$, fix $x\in 2^\om$ and $n,m<\om$.  The following 
set is dense:
$$D_{x,n,m}=\{p\in\poset(A): \exists k\;\; \lq[x\res k]\sq U_{nm}\rq\in p\rmor
\lq x\notin U_{nm}\rq\in p\}$$
To see this note that if $\lq x\notin U_{nm}\rq$ is not in $p$ we can 
always find
$k$ large enough so that  $p\cup \{\lq [x\res k]\sq U_{nm}\rq\}$ is
a consistent set of sentences.
Now suppose $x\in U^G_{nm}$, then for some $k$ we have
that $\lq[x\res k]\sq U_{nm}\rq\in G$ and hence by consistency,
$\lq x\notin U_{nm}\rq\notin G$.  On the otherhand, if 
$\lq x\notin U_{nm}\rq\notin G$, then since $D_{x,n,m}$ is dense
for some $k$ we have that $\lq[x\res k]\sq U_{nm}\rq\in G$ and hence
$x\in U^G_{nm}$.

To prove $(2)$ note that the following set is dense:
$$D_{x,n}=\{p\in\poset(A)\;:\; \exists k\;\; \lq x\notin U_{nk}\rq \in p\rmor
\lq x\in\cap_{m<\om}U_{nm}\rq\in p\}$$
To see this note that if $\lq x\in\cap_{m<\om}U_{nm}\rq\notin p$, then
for large $k$ (so that $U_{nk}$ is not mentioned in $p$), the sentences
$p\cup \{\lq x\notin U_{nk}\rq\}$ are consistent. 

To prove $(3)$ note that if $x\in A$ then the following is dense:
$$D_{x}=\{p\in\poset(A)\;:\; \exists n\; \lq x\in\cap_{m<\om}U_{nm}\rq\in p\}$$
and we can only assert $\lq x\in\cap_{m<\om}U_{nm}\rq$ for $x\in A$. 
\qed

Note that it follows from the Lemma that $A\cap V=(\cup_{n<\om}W_n^G)\cap V$
and so that $A$ is a $\Sigthree$ relative to the ground model reals.

\begin{lemma}\label{lemccc} $\poset(A)$ is ccc.
\end{lemma}
\proof
This is a standard $\Delta$-systems argument.  Suppose two conditions
$p$ and $q$ agree on all sentences of the form:
 $$\lq [s]\sq U_{nm}\rq $$
and also they agree on all sentences of the form:
 $$\lq x\in \cap_{m<\om}U_{nm}\rq \rmor \lq x\notin U_{nm}\rq $$
whenever $x$ is mentioned in both $p$ and $q$.  Then $p\cup q$ is
consistent.
\qed

\bigskip
Next we must prove that $\poset(A)$ does not add a dominating real.

Working in $V$, for $Y\sq 2^\om$ countable define
$p\in \poset(A)_Y$ iff $p\in\poset(A)$ and 
$$\forall x,n,k\;\; 
(\lq x\notin U_{nk}\rq \in p\rmor 
 \lq x\in\cap_{m<\om} U_{nm}\rq \in p)\to x\in Y\}.$$
Or in otherwords, $\poset(A)_Y$ are the conditions in $\poset(A)$ which
only mention elements of $Y$.   

\begin{lemma}\label{easy}
Suppose $p\in\poset(A)$ and $q\in\poset(A)_Y$.  Then 

$p$ and $q$ are compatible iff $r$ and $q$ are compatible

where 
$$r=p\setminus\{\lq x\in\cap_{m<\om}U_{nm}\rq \;:\; x\notin Y, n<\om\}$$
\end{lemma}
\proof
Incompatibility cannot arise between sentences of type (1) and (3).  That
is, any 
pair of the form:
$$\lq [s]\sq U_{nm}\rq,\;\; \lq x\in\cap_{m<\om}U_{nm}\rq$$
is consistent.
It follows that the $\lq x\in\cap_{m<\om}U_{nm}\rq\in p$ for
which  $x\notin Y$ cannot conflict with the sentences of $q$ since
by definition $q$ cannot mention any $x$ which is not in $Y$.
\qed

\bigskip\noindent Define. $T=(p,(t_i,n_i,m_i:i<N))$ is a $Y$-template iff
\begin{enumerate}
 \item $p\in \poset(A)_Y$, $t_i\in 2^{<\om}$, $n_i,m_i,N\in \om$, 
 \item if $\lq y\in \cap_{m<\om} U_{n_im}\rq \in p$, then $y\notin [t_i]$, and
 \item if $\lq [s]\sq U_{n_im_i}\rq\in p$, then $[s]\cap [t_i]=\emptyset$.
\end{enumerate}

\bigskip\noindent Define. For $\vec{x}=(x_i:i<N)\in \prod_{i<N}[t_i]$
$$p(\vec{x})=p\cup \{\lq x_i\notin U_{n_im_i}\rq\;:\; i<N\}$$

\bigskip
Note that by the definition of $Y$-template that $p(\vec{x})\in\poset(A)$,
i.e., is consistent, for every $\vec{x}\in \prod_{i<N}[t_i]$.

\begin{lemma}\label{mainlem}
Suppose that $\forces \tau\in\om$,
there exists $\Sigma\sq \poset(A)_Y$ a maximal antichain deciding $\tau$,
and $(p,(t_i,n_i,m_i:i<N))$ is a $Y$-template. Then there exists $k<\om$
so that for every $\vec{x}\in \prod_{i<N}[t_i]$ there exists
$q\in \poset(A)_Y$ such that $p(\vec{x})\cup q\in\poset(A)$ and 
$q\forces \tau<k$.
\end{lemma}
\proof
For $q\in\poset(A)_Y$ define
$$U_q=\{\vec{x}\in\prod_{i<N}[t_i]\;:\; p(\vec{x})\cup q\in\poset(A)\}$$
Note that $U_q$ is open.  To see this, suppose $\vec{x}\in U_q$ so
that $p(\vec{x})\cup q\in\poset(A)$.  Note that although some $x_i$ might be
in $Y$ it can't be that $\lq x_i\notin U_{n_im_i}\rq\in p(\vec{x})$ and
$\lq x_i\in\cap_{m<\om}U_{n_im}\rq\in q$, because they are compatible. 
Hence, there must be a sufficiently small neighborhood of $x_i$ say
$t_i'=x_i\res k_i\supseteq t_i$ with the properties that  
\begin{enumerate}
\item if $\lq z\in \cap_{m<\om} U_{n_im}\rq \in p\cup q$, then $z\notin
[t_i']$, and
\item if $\lq [s]\sq U_{n_im_i}\rq\in p\cup q$, then $[s]\cap
[t_i']=\emptyset$.
\end{enumerate}
Hence, $\vec{x}\in\prod_{i<N}[t_i']\sq U_q$.

Now since  $\Sigma\sq \poset(A)_Y$ is a maximal
antichain we know that 
$$\cup\{U_q:q\in\Sigma\}=\prod_{i<N}[t_i] $$
So by compactness since each $U_q$ is open, 
there exists a finite $F\sq\Sigma$ such that
$$\cup\{U_q:q\in F\}= \prod_{i<N}[t_i]$$
and since each $q\in\Sigma$ decides $\tau$, the Lemma follows.
\qed

In order to prove the full result we must show that the iteration
does not add a dominating real.  To do this we prove the following stronger
property (see Bartoszynski and Judah \cite{barto} definition 6.4.4):

\begin{lemma}
The poset $\poset(A)$ is really $\sqbdd$-good, i.e., for every name
$\tau$ for an element of $\om^\om$ there exists $g\in\om^\om$ such 
that for any $x\in\om^\om$ if
there exists $p\in\poset(A)$ such that 
$p\forces\forall^\infty n\; x(n)<\tau(n)$,
then $\forall^\infty n\; x(n)<g(n)$. 
\end{lemma}
\proof

Suppose that $\forces \tau\in\om^\om$.  Let $Y\sq 2^\om$ be countable so
that for every $n<\om$ there exists a maximal antichain $\Sigma\sq \poset(A)_Y$
which decides $\tau(n)$. List all 
$Y$-templates as $(T_n:n<\om)$.  
By Lemma \ref{mainlem}
there exists $g\in\om^\om$ with the property that for every
$l<\om$ and $n<l$ if 
$$T_n=(p,(t_i,n_i,m_i:i<N))$$ 
then for every 
$\vec{x}\in \prod_{i<N}[t_i]$ there exists $q\in \poset(A)_Y$ such
that $p(\vec{x})\cup q\in\poset(A)$ and $q\forces \tau(l)<g(l)$.
(To get $g(l)$ apply Lemma \ref{mainlem} to $\tau=\tau(l)$ and
each of the templates $(T_n:n<l)$ and then take $g(l)$ to be the maximum
of all the $k's$.)
 
Now suppose that $p_0\forces \forall l>l_0 \;\;x(l)<\tau(l)$ and
$$p_0=p\cup\{z_i\in \cap_{m<\om}U_{n'_i,m}: i< N'\}\cup 
\{x_i\notin U_{n_im_i}\;:\; i< N\}$$
where $p\in \poset(A)_Y$ and $z_i,x_i\notin Y$.  
 
Take $t_i$ sufficiently long so that $t_i\sq x_i$ and  
$$T=(p, (t_i,n_i,m_i:i<N))$$ 
is a $Y$-template.  Assume that $l_0$ is sufficiently large so
that $T=T_k$ for some $k<l_0$. 
By our construction for each $l>l_0$, there exists $q\in\poset(A)_Y$
such that $p(\vec{x})\cup q\in\poset(A)$ and $q\forces \tau(l)<g(l)$. But
by Lemma \ref{easy} this means that $p_0\cup q\in \poset(A)$
and hence $x(l)<g(l)$. 
 
\qed

The above proof is similar to that of Lemma 6.5.8 \cite{barto}.

Now we prove Theorem \ref{zapletal}.  Starting with a model of CH we
iterate with finite support $\om_2$ times 
$$\poset_{\alpha+1}=\poset_\alpha*\name{\poset}(\name{A_\alpha})$$
where we dovetail to list all $A\sq 2^\om$ of size $\om_1$ in the
final model.  Since the finite support iteration of really 
$\sqbdd$-good ccc forcing adds no dominating real
(see Bartoszynski and Judah \cite{barto} Theorem 6.5.4), 
we have that in the resulting model that $\bb=\om_1$.
On the other hand by Lemma \ref{one} we have that
$\zz=\om_2$. 
\qed

\bigskip
Define (see Zapletal \cite{zap} Appendix A)
$$\sn=\min\{|X|: X\sq \trees, 
\forall  A \;\;\Sig\; \;\;X\cap A\not= X\cap WF\}$$
where $\trees$ is the set of $\om$-trees and $WF$ is the set of
well-founded trees.  An equivalent definition is:
$$\sn=\min\{|X| \;:\; X\sq 2^\om\; \exists A\;\;\Sig\;\;
\forall B\;\;\Pig\;\; X\cap A\not= X\cap B\}$$ 
    
The equivalence is easy to show because 
the set of well-founded trees is a universal $\Pig$ set.
It is not hard to see that $\zz\leq\sn$. 
So we have the relative consistency of $\bb<\sn$.  

The following proposition is mostly due to Rothberger \cite{roth}. It implies
that we must go up to at least the third level of the Borel hierarchy to get 
the consistency of $\bb<\sn$.

\begin{prop} For $\kappa$ an infinite cardinal the following are equivalent:

\begin{enumerate}

\item $\bb>\kappa$

\item For all $X\sq 2^\om$ with $|X|\leq\kappa$ and for all 
$\Sig$ sets $A\sq 2^\om$ there exists a
$\Sigtwo$ set $B\sq 2^\om$ such
that $X\cap A=X\cap B$.

\item For all $X\sq 2^\om$ with $|X|\leq\kappa$ and for all 
$\Sigtwo$ sets $A\sq 2^\om$ there exists a
$\Pigtwo$ set $B\sq 2^\om$ such
that $X\cap A=X\cap B$.

\item For all $X\sq 2^\om$ with $|X|\leq\kappa$ and for all 
countable  $A\sq X$ there exists a
$\Pigtwo$ set $B\sq 2^\om$ such
that $A=X\cap B$.

\end{enumerate}
\end{prop}

\proof

$(2)\to (3)$ and $(3)\to (4)$ are trivial. 

To see $(1)\to (2)$ let
$$A=\{x\in 2^\om\;:\; \exists y\in \om^\om\;\;(x,y)\in C\}$$
where $C\sq 2^\om\times \om^\om$ is closed. 
Suppose that $A\cap X=\{x_\alpha:\alpha<\kappa\}$.  Choose
$y_\alpha\in \om^\om$ so that $(x_\alpha,y_\alpha)\in C$ for
each $\alpha<\kappa$. Since $\bb>\kappa$ we can choose $z_n\in\om^\om$
for $n<\om$ so that for all $\alpha<\kappa$ there exists $n<\om$ 
with $y_\alpha\leq z_n$ (pointwise).  Define
$$C_n=\{(x,y)\in C: y\leq z_n\}$$
$C_n$ is compact and therefore so is its projection:
$$A_n=\{x\in 2^\om \;:\;\exists y \;(x,y)\in C_n\}$$
But $A\cap X=\cup_{n<\om}A_n\cap X$.

To see $(4)\to (1)$ let $X\sq \om^\om$
with $|X|=\kappa$.  Now since $\om^\om$ is homeomorphic to $\infsub$
and $\infsub\sq P(\om)\simeq 2^\om$ by applying $(4)$ we can
find a $\Pigtwo$ set $G\sq P(\om)$ such that
$$G\cap (X\cup \finsub)=\finsub$$
But note that $F=P(\om)\setminus G$ is a $\sigma$-compact set which is
disjoint from $\finsub$, i.e. a subset of $\infsub\simeq \om^\om$
and covers $X$.  
But is easy to show that for any $\sigma$-compact subset $F$ of $\om^\om$
there exists $f\in \om^\om$ such that $g\leq^*f$ for all $g\in F$.
\qed

\bigskip
Remark.
One way to get the consistency of $\bb<\zz<\sn$ is as follows:
Start with a ground model of $2^\om=\om_1$, $2^{\om_1}=\om_2$, and
$2^{\om_2}=\om_{17}$.  Do a finite support iteration of
$\poset(A_\alpha)$ for $\alpha<\om_3$, so
that for each $\alpha$ either $A_{\alpha}=A$ the universal $\Sig$-set or
$|A_\alpha|=\om_1$ as in the above proof. In the final
model we will have $\bb=\om_1$ since it is an iteration of
really $\sqbdd$-good ccc partial orders.  Also we will have
$\zz\leq \om_2$ because $2^{\om_2}=\om_{17}$ and $2^\om=\om_3$.
We also have $\zz\geq\om_2$ because of dovetailing over all $|A|=\om_1$.
And we will have $\sn=\om_3=\cc$ because we have cofinally used the
universal $\Sig$-set.

\bigskip
\bigskip
The following Theorem answers a question of Dan Mauldin
(see \cite{prob} problem 7.8).

\begin{theorem}\label{mauldin}
It is relatively consistent with ZFC that there exist a separable metric
space $X$ such that the Borel order of $X$ is bounded, but not every
relatively analytic subset of $X$ is Borel in $X$.
\end{theorem}

\proof
We use almost exactly the same partial order but with one crucial difference.
Instead of using arbitrary subsets $A\sq 2^\om$ we let 
$B\sq 2^\om$ be a fixed universal $\Pigthree$ set. The partial order
$\pos$ is Borel, ccc, and adds a generic $\Sigthree$ set whose
intersection with the ground model is the same as $B$'s with the ground 
model. 

\bigskip\noindent
Define. A partially ordered set $\poset$ is very Souslin iff
\begin{enumerate}

\item $\poset$ is ccc,

\item $\poset$, $\;\leq$, $\;\{(p,q)\in \poset^2: p, q \mbox{ incompatible }
\}$
are $\Sig$, and

\item $\{\Sigma\in \poset^\om \;:\; \Sigma 
\mbox{ enumerates a maximal antichain }\}$
is $\Sig$.

\end{enumerate}

\bigskip
We will need the following Lemma:

\begin{lemma} (Zapletal \cite{zap} see Appendix C, Lemmas
C.0.14 and C.0.17)\label{zaptwo}
Suppose $\poset$ is a very Souslin real partial order and
$\poset^{\om_2}$ the countable support iteration of $\poset$. Then
$$V^{\poset^{\om_2}}\models \sn =\om_1$$
\end{lemma}

Clearly this means that partial order $\poset(A)$  is not very Souslin even
when $A$ is taken to be analytic (so it is Souslin). However if we change $A$
to make it Borel, then it is very Souslin.

\begin{lemma}
The partial order $\pos$ is very Souslin.
\end{lemma}
\proof
The following sets are Borel:
\begin{enumerate}
 \item $\pos$
 \item $\{(p,q)\in \pos\times\pos\;:\; p\sq q\}$
 \item $\{(p,q)\in \pos\times\pos\;:\; \mbox{ $p$ and  $q$ are incompatible }
   \}$
 \item $\{(p,Y)\;:\; Y\in [2^\om]^\om \rmand p\in \pos_Y\}$ 
 \item $\{((T_n:n<\om),Y)\;:\; Y\in [2^\om]^\om $ and 
       $\{T_n:n<\om\}=$ all $Y$-templates $\}$        
\end{enumerate}

\bigskip
\noindent Next we verify that being a maximal antichain in $\pos$ is $\Sig$. 

\bigskip\noindent{\bf Claim}.
$\Sigma\sq\pos$ is a maximal antichain iff
\begin{enumerate}
\item $\Sigma$ is an antichain and
\item there exists $Y\sq 2^\om$ countable and $(T_n:n<\om)$
such that 
  \begin{itemize}
       \item $\Sigma\sq \pos_Y$ and 
       \item $(T_n:n<\om)$ enumerates the set of all  $Y$-templates 
  \end{itemize}
and for all $n$ if $T_n=(p,(t_i,n_i,m_i:i<N))$, then there exists 
\par $K$, $\;(t_i^j:j<K)$, and $(q_j:j<K)$ such that
  \begin{enumerate}
    \item $\prod_{i<N}[t_i]=\cup_{j<K}\prod_{i<N} [t_i^j]$
    \item $q_j\in \Sigma$     
    \item $q_j\cup p \in \pos$ 
    \item $\lq y\in\cap_{m<\om}U_{n_i,m}\rq \in q_j \to y\notin [t_i^j]$
    \item $\lq [s]\sq U_{n_im_i}\rq\in q_j \to [t_i^j]\cap [s]=\emptyset$ 
  \end{enumerate} 
\end{enumerate}
\proof
Condition (2) is just a detailed restatement of Lemma \ref{mainlem} and 
its proof. It guarantees by Lemma \ref{easy} that every $p\in\pos$ is 
compatible with some $q\in \Sigma$. 

This proves the claim and the lemma easily follows.
\qed

Hence by Zapletal's Lemma \ref{zaptwo} if we iterated $\pos$ with countable
support $\om_2$ times then in the resulting model $\sn=\om_1$. 
Hence there is some
$X\sq 2^\om$ of size $\om_1$ with a relatively analytic set which is not
relatively coanalytic.  (Actually the proof of Lemma \ref{zaptwo} 
shows that the ground model
reals would do for such an $X$).
But note that every $\Pigthree$ set occurs as a cross section of our
universal $\Pigthree$-set $B$ and by
Lemma \ref{one} becomes $\Sigthree$ with respect
to the ground model.  Hence it is easy to see 
that for every $X\sq 2^\om$ of size $\om_1$ for
every $\Sigthree$ $B$ there exists a $\Pigthree$ $C$ such that
$X\cap B=X\cap C$.  This proves Theorem \ref{mauldin}.
\qed

\begin{flushleft}
Arnold W. Miller \\
miller@math.wisc.edu \\
http://www.math.wisc.edu/$\sim$miller\\
University of Wisconsin-Madison \\
Department of Mathematics, Van Vleck Hall \\
480 Lincoln Drive \\
Madison, Wisconsin 53706-1388 \\
\end{flushleft}

\appendix

\newpage
\begin{center}
Appendix
\end{center}

(Not intended for publication, electronic version only.)

Our first proof of $\bb<\sn$ used large cardinals
and the following Lemma:

\begin{lemma}\label{zaplem}
(Zapletal \cite{zap} Thm 5.4.12) (LC)
Suppose $\poset$ is a real, proper, universally
Baire forcing such that 
$$V^\poset\models V\cap \om^\om\mbox{is unbounded in }
\leq^*$$
Then 
$$V^{\poset^{\om_2}}\models V\cap \om^\om\mbox{is unbounded in }
\leq^*$$
where $\poset^{\om_2}$ stands for the $\om_2$ iteration with
countable support of $\poset$.
\end{lemma}
The hypothesis (LC) stands for large cardinals, for example, unboundedly
many measurable Woodin cardinals would be enough.
In otherwords for a nice enough forcing, not adding a dominating real is
preserved by the iteration.  It is easy to get a two step iteration
so that neither step adds a dominating real but the two steps do. 
For example, force $\om_1$-Cohen reals followed by the Heckler
partial order of the ground model.

Fix $A\sq 2^\om$ a universal $\Sigma^1_1$ set, i.e., it
is lightface $\Sigma_1^1$ and every boldface $\Sig$ occurs as
a cross section via some effective homeomorphism of $2^\om\times 2^\om$
and  $2^\om$.
In this case the partial order 
$\poset(A)$ is $\Sig$, ccc, and determined
by a real - so it satisfies the hypothesis of the Lemma.

\end{document}